\documentclass[12pt,a4paper]{article}

\usepackage{a4}

\usepackage[intlimits,tbtags]{amsmath} 
\usepackage{amsfonts,amsthm}
\usepackage{amscd}

\usepackage{cite,latexsym} 

\newcommand{\E}{\mathrm{e}} 
\newcommand{\I}{\mathrm{i}}
 
\newcommand{\tr}{\triangleright}

\newcommand{\otimeshat}{\operatorname{\hat{\otimes}}}

\newcommand{\Xcal}{\mathcal{X}}

\newcommand{\Ocal}{\mathcal{O}}

\newcommand{\Ket}[1]{\lvert #1\rangle} 
\newcommand{\Bra}[1]{\langle #1\rvert}

\newtheorem{Theorem}{Theorem}

\newcommand{\h}{\hbar}

\newcommand{\Uhsl}{{\mathcal{U}_\h(\mathrm{sl}_2)}}
\newcommand{\Usl}{{\mathcal{U}(\mathrm{sl}_2)}}

\newcommand{\Uhsu}{{\mathcal{U}_\h(\mathrm{su}_2)}}
\newcommand{\Uhg}{{\mathcal{U}_\h(\mathfrak{g})}}

\newcommand{\Ug}{{\mathcal{U}(\mathfrak{g})}}

\newcommand{\Uhso}{{\mathcal{U}_{\h}(\mathrm{so}_4)}}
\newcommand{\Uso}{{\mathcal{U}(\mathrm{so}_4)}}

\newcommand{\CGs}[6]{
  \bigl(\begin{smallmatrix}#1\! &  #2\vphantom{#3} \\
    #4\! & #5\vphantom{#6} \end{smallmatrix} \!\!\bigm|\!\!
  \begin{smallmatrix}
    \vphantom{#1#2}#3 \\ \vphantom{#4#5}#6
  \end{smallmatrix}\bigr)}

\newcommand{\CGqs}[6]{\CGs{#1}{#2}{#3}{#4}{#5}{#6}_{\! q}}

\newcommand{\qBinom}[2]{{\genfrac{[}{]}{0pt}{1}{#1}{#2}}_{q^{-2}}}

\begin{document}

\rightline{IUB-TP/2003-10}
%\rightline{MPI-ThP/2002-43}
\vspace{3em}
\begin{center}

  {\Large{\bf Perturbative Symmetries on Noncommutative Spaces}}

\vspace{3em}
\textbf{Christian Blohmann}
\\[1em]
International University Bremen, School of Engineering and Science\\
Campus Ring 1, 28759 Bremen, Germany
\end{center}

\vspace{1em}

\begin{abstract}
  Perturbative deformations of symmetry structures on noncommutative
  spaces are studied in view of noncommutative quantum field theories.
  The rigidity of enveloping algebras of semi-simple Lie algebras with
  respect to formal deformations is reviewed in the context of star
  products. It is shown that rigidity of symmetry algebras extends to
  rigidity of the action of the symmetry on the space.  This implies
  that the noncommutative spaces considered can be realized as star
  products by particular ordering prescriptions which are compatible
  with the symmetry. These symmetry preserving ordering prescriptions
  are calculated for the quantum plane and four-dimensional quantum
  Euclidean space. Using these ordering prescriptions greatly
  facilitates the construction of invariant Lagrangians for quantum
  field theory on noncommutative spaces with a deformed symmetry.
%  \keywords{noncommutative geometry, deformation theory, star
%    products, quantum spaces, quantum symmetries}
\end{abstract}

\section{Introduction}

If the geometry of a physical space is noncommutative at energies
accessible by current accelerators \cite{Seiberg:1999} the
noncommutativity can only be small.  This suggests to describe
noncommutative spaces as perturbative deformations of ordinary,
commutative spaces. If such a small deformation is to have
controllably small effects, it must depend in some sense smoothly on a
deformation parameter.  Ideally, the deformation of a physical
quantity which results from a perturbation of the geometry is given by
a convergent perturbation series in powers of the deformation
parameter. But even if convergence in this strong sense is rigorously
not possible, as for quantum field theory or deformation quantization,
the perturbation series may still be useful in an algorithmic or
algebraic sense.

The algebraic aspects of a deformation can be separated from the
analytic questions of continuity and convergence by considering formal
power series. In such a framework, a noncommutative space can be
described as formal deformation of the algebra of functions on the
space manifold in the sense of Gerstenhaber \cite{Gerstenhaber:1964}:
The deformed algebra is an algebra over the ring of formal power
series in the perturbation parameter, which is in zeroth order
isomorphic to the undeformed algebra. As it turns out, formal power
series are the natural setting for the construction of gauge theories
on general noncommutative geometries, which have found a solid
formulation \cite{Madore:2000b,Jurco:2001} within the framework of
deformation quantization. For reviews of noncommutative field theories
see \cite{Douglas:2001} and \cite{Szabo:2001}.

The perturbative approach to noncommutative field theories, that is,
expanding the product of noncommutative quantum fields in the
noncommutativity parameters and relating the noncommutative gauge
potentials and fields to their ordinary, commutative counterparts via
the Seiberg--Witten map has put quantum theories on noncommutative
spaces within the range of phenomenological considerations: A minimal
noncommutative extension of the standard model was
formulated \cite{Calmet:2001} the effects of noncommutative geometry
on magnetic and electric moments were
studied \cite{Iltan:2003,Minkowski:2003} noncommutative
neutrino-photon coupling with possible astrophysical implications was
investigated \cite{Schupp:2002} the OPAL collaboration has started
looking for noncommutative signatures in electron positron pair
annihilation \cite{Abbiendi:2003} just to name some recent examples.
For a review on the phenomenological implications of noncommutative
geometry see \cite{Hinchliffe:2002}.

As important as the spaces are the symmetries which act on them.
Deforming a space algebra which transforms covariantly under a
symmetry Lie group will in general break the symmetry.  For example,
the noncommutative geometry which most papers on noncommutative field
theory have considered, where the commutator of the space-time
observables $[X_\mu, X_\nu] = \theta_{\mu\nu}$ is a constant
antisymmetric matrix, breaks Lorentz symmetry. Physically, this has to
be expected as the constant commutator can be viewed as due to a
constant background field, in string theory a constant $B$-field on a
D-brane.  It could be argued that if the noncommutativity parameters
$\theta_{\mu\nu}$ are small, the violation of Lorentz symmetry is only
small, too. However, on the level of regularization of loop diagrams
the noncommutativity leads to an interdependence of ultra-violet and
infra-red cutoff scales \cite{Minwalla:1999,Matusis:2000}. This UV/IR
mixing seems to put even large scale Lorentz symmetry and weakened
notions of locality of noncommutative quantum field theory into
doubt \cite{Alvarez-Gaume:2003}. But as yet, UV/IR mixing was
investigated in detail only for the case of constant
$\theta_{\mu\nu}$.

The appearance of UV/IR mixing seems to indicate, that the breaking of
symmetries which happens when a space is noncommutatively deformed
with constant $\theta_{\mu\nu}$ is not under good control. In a
self-contained theory it would be reasonable to expect
$\theta_{\mu\nu}$ to become itself a field, which transforms
covariantly with respect to a perturbative deformation of space-time
symmetry.  The existence of a deformed symmetry structure would be a
big advantage for phenomenological considerations. It would allow to
include in the perturbative approach the changes induced by
noncommutativity to those physical concepts that are tied to
space-time symmetry, such as energy-momentum conservation, Lorentz
invariance, independence of in and out states etc.

Another motivation to study noncommutative spaces with deformed
symmetries has emerged recently from the attempts to explain the
observation \cite{Bird:1995} of cosmic rays of energy beyond the
spectral cutoff (the Greisen--Zatsepin--Kuzmin limit) which is
expected due to interaction with the cosmic microwave background. The
often proposed explanation of such ultra high energy rays by vacuum
dispersion relations, that is, the dependence of the speed of light on
the wavelength, was shown by Amelino-Camelia to be reconcilable in
principle with the observer independence of the laws of physics
\cite{Amelino-Camelia:2000}. This leads to a deformation of special
relativity by the assumption that there is not only an observer
invariant velocity but also an observer invariant length, the Planck
length, which plays the role of the deformation parameter. This
proposition, now called doubly special relativity, has initiated a
large number of active studies from both, the mathematical and the
phenomenological viewpoint. (For an overview see
\cite{Amelino-Camelia:2002}.) Realizations of doubly special
relativity can lead to noncommutative deformations of the space-time
with a deformed Lorentz symmetry \cite{Agostini:2003}.

The purpose of this paper is to study formal perturbative deformations
of symmetry structures on noncommutative spaces. On ordinary
commutative spaces symmetry structures can be described by Lie
algebras or, equivalently, their enveloping algebras.  If the Lie
algebra is semi-simple as for most interesting cases in physics, the
relation of an enveloping algebra to its deformation turns out to be
surprisingly simple: The two algebras are isomorphic.  More precisely,
it can be shown by homological arguments \cite{Gerstenhaber:1964} that
the enveloping algebra of a semi-simple Lie algebra is \emph{rigid},
that is, over formal power series \emph{any} deformation is isomorphic
to the undeformed algebra.  We recall this result in
theorem~\ref{th:RigidAlgebra}.

In addition to the symmetry algebra we need the action of the symmetry
algebra on the space in order to describe the symmetry structure
completely. \emph{A priori}, even though the symmetry algebra is
rigid, the action could be truly deformed. However, using the same
homological methods as before, we show, that the action is rigid, as
well: Over formal power series, the space with an arbitrarily deformed
action is isomorphic as module to the space with the usual action of
the enveloping Lie algebra by differential operators. This result is
stated in theorem~\ref{th:RigidAction}. In the context of star
products vector space isomorphisms between the deformed and the
undeformed space are often referred to as ordering prescriptions. In
this language theorem~\ref{th:RigidAction} shows that there are
particular ordering prescriptions which are compatible with the
symmetry structure of the space.

Since the widely studied case of constant $\theta_{\mu\nu}$ does not
allow for a perturbative deformation of Lorentz symmetry, it cannot
serve as example for these general theorems. The standard examples for
noncommutative spaces with deformed symmetry structures are quantum
spaces \cite{Manin:1988,Faddeev:1990,Carow-Watamura:1990}. They carry a
covariant representation of the Drinfeld--Jimbo
deformation \cite{Drinfeld:1985,Jimbo:1985} of the enveloping symmetry
algebra. In fact, rigidity theorem~\ref{th:RigidAlgebra} which we
review here lies at the core of Drinfeld's pioneering work on quantum
enveloping algebras \cite{Drinfeld:1989,Drinfeld:1989b}.

While from a mathematical point of view isomorphic objects are often
identified, an isomorphism can change the physical interpretation of
the symmetry structure.  Finding the isomorphisms between the deformed
and undeformed structures on an algebraic level, however, is a
difficult computational problem because the homological methods which
are used to prove the rigidity theorems, although elegant, are
inherently non-constructive. We will use representation theory to
reduce the algebraic problem to matrix calculations.  This approach
works well for cases, where the representation theory is well
understood such as for quantum spaces and quantum algebras.

We will take the quantum plane with its $\Uhsu$-symmetry as guiding
example to demonstrate these representation theoretic methods. The
construction of isomorphisms between enveloping algebras and their
quantum deformations is text-book material (e.g. Sec.~6.1.3 of
\cite{Schmuedgen}). The isomorphism of the module structures has
received less attention, but it is of importance for the realization
of noncommutative spaces by star products: As in deformation
quantization, the multiplication map of a given noncommutative space
algebra is often transferred to a commutative function algebra using
an ordering prescription, by which the spaces are identified as vector
spaces. The rigidity of the module structure as formulated in
theorem~\ref{th:RigidAction} tells us that there is an ordering
prescription which is not only an isomorphism of vector spaces but
also of modules. The fact that the deformed and undeformed spaces are
isomorphic as modules will only be obscured by most ordering
prescriptions, such as the popular normal ordering and the symmetric
ordering.

The main result of this paper is the calculation of the symmetry
preserving ordering prescription for quantum Euclidean four-space in
Eq.~\eqref{eq:phi3}. The result is expressed in terms of the deformed
and undeformed binomial and Clebsch--Gordan coefficients, that is, in
terms of basic hypergeometric series. In this sense the representation
theoretic approach profits from the computational effort that has gone
into the calculation of the $q$-Clebsch--Gordan coefficients. Trying
to redo the calculation which leads to Eq.~\eqref{eq:phi3} in a
recursive fashion order by order in the deformation parameter, one
would quickly learn that $q$-hypergeometric functions are an extremely
efficient way to describe the complex combinatorics of partitions, the
reason for which they were first introduced by Euler in the 18th
century.

For a rigorous outline of the basics of $\h$-adic topologies on
algebras in the context of deformation theory, which are used in this
article, we refer the reader to \cite{Kassel}, Ch.~XVI.

\section{Space Deformations and Symmetries}
\label{sec:Actions}

\subsection{A guiding example}

As guiding example for a noncommutatively deformed space with a
symmetry structure let us take the $\h$-adic quantum plane, the
$\h$-adic complex algebra with two generators $\hat{x}$, $\hat{y}$ and
commutation relation $\hat{x}\hat{y} = \E^\h \hat{y}\hat{x}$.
Expanding this relation in orders of the deformation parameter, we see
that in zeroth order the generators $\hat{x}$ and $\hat{y}$ commute.
Hence, the polynomial algebra in two variables and the quantum plane,
\begin{equation}
\label{eq:PlanesDef}
 \Xcal := \mathbb{C}[ x, y ]
 \quad\text{and}\quad
 \Xcal_\h := \mathbb{C}\langle\hat{x},\hat{y}\rangle[[\h]] /
      \langle \hat{x}\hat{y} = \E^\h \hat{y}\hat{x} \rangle \,,  
\end{equation}
are isomorphic modulo $\h$. That is, there is an isomorphism of
algebras $\xi: \Xcal \rightarrow \Xcal_\h/\h\Xcal_\h$ which is defined
on the generators as $\xi(x) = \hat{x}$, $\xi(y) = \hat{y}$.

The isomorphism of algebras $\xi$ is \emph{a fortiori} an isomorphism
of vector spaces and can be extended to formal power series yielding a
$\mathbb{C}[[\h]]$-linear isomorphism of $\h$-adic vector spaces
$\varphi: \Xcal[[\h]] \rightarrow \Xcal_\h$. Such an extension
$\varphi$, which we will call an ordering prescription, is not unique.
For example, the image of the quadratic term $xy$ could be equally
defined as $\varphi(xy) := \hat{x}\hat{y}$ (normal ordering) or as
$\varphi(xy) := \tfrac{1}{2}(\hat{x}\hat{y} + \hat{y} \hat{x} )$
(symmetric ordering), just to take two popular ordering prescriptions.

Using the ordering prescription, we can transfer the noncommutative
multiplication map $\mu_\h$ of $\Xcal_\h$ to $\Xcal[[\h]]$ by
requiring
\begin{equation}
\label{eq:StarProd}
\begin{CD}
  \Xcal[[\h]] \otimeshat \Xcal[[\h]] 
  @>{\varphi \otimes \varphi }>> \Xcal_\h \otimeshat \Xcal_\h \\
  @VV{\mu_\varphi}V @VV{\mu_\h}V \\ 
  \Xcal[[\h]] @>{\varphi }>> \Xcal_\h
\end{CD}
\end{equation}
to be a commutative diagram, where $\otimeshat$ denotes the
topological tensor product. We will call the transferred
multiplication map
\begin{equation}
\label{eq:starmult}
  \mu_\varphi :=
  \varphi^{-1} \circ \mu_\h \circ (\varphi \otimes \varphi )
\end{equation}
a star product, denoted by $\mu_\varphi(x \otimes x') \equiv x \star
x'$. By construction, the vector space $\Xcal[[\h]]$ equipped with
this star product is now isomorphic as algebra to $\Xcal_\h$,
$(\Xcal[[\h]], \mu_\varphi) \cong (\Xcal_\h, \mu_\h)$. While a
different ordering prescription $\varphi'$ will in general yield a
different multiplication map $\mu_{\varphi'} \neq \mu_\varphi$, the
algebras will be isomorphic, $(\Xcal[[\h]], \mu_{\varphi'}) \cong
(\Xcal[[\h]], \mu_\varphi)$, with $\varphi^{-1} \circ \varphi'$ being
an isomorphism.

At this point the construction of the star product seems somewhat
vacuous. The reason for transferring the noncommutativity to the
ordinary function space is the additional structure on $\Xcal$ which
might not be present on $\Xcal_\h$, such as a differential calculus,
integration, or a symmetry.

The function algebra of the plane can be equipped with a $\Usl$
symmetry structure. The Lie algebra $\mathrm{sl}_2$ of the special
linear group in two dimensions is generated by the Cartan--Weyl
generators $E$, $H$, and $F$ with relations $[H,E]=2E$, $[H,F]=-2F$,
$[E,F] = H$. By definition, the elements of a Lie algebra act on the
function space $\Xcal$ as derivations, so we can represent the
generators by first order differential operators
\begin{equation}
  E = y\partial_x \,,\quad H = y\partial_y - x \partial_x
  \,,\quad F = x \partial_y \,. 
\end{equation}
This symmetry structure can be deformed together with the deformation
of the algebra of the plane into the quantum plane. The deformed
symmetry algebra $\Uhsl$ is generated by $\hat{E}$, $\hat{H}$, and
$\hat{F}$ with relations $[\hat{H},\hat{E}]=2\hat{E}$,
$[\hat{H},\hat{F}]=-2\hat{F}$, and 
\begin{equation}
\label{eq:qCommutation}
  [\hat{E},\hat{F}]
  = \frac{\E^{\h \hat{H}}-\E^{-\h \hat{H}}}{\E^\h-\E^{-\h}}\,.
\end{equation}
We can define the action of $\hat{E}$, $\hat{H}$, and $\hat{F}$ on the
generators $\hat{x}$ and $\hat{y}$ of the quantum plane exactly as for
the undeformed case: $\hat{H}$ is diagonal, $\hat{H}\tr x = -x$,
$\hat{H}\tr y = y$, and $\hat{E}$, $\hat{F}$ are ladder operators,
$\hat{E} \tr x = y$, $\hat{E} \tr y = 0$, $\hat{F} \tr y = x$,
$\hat{F} \tr x = 0$. However, while $\hat{H}$ still acts on $\Xcal_h$
as derivation, we have to modify the Leibniz rule for $\hat{E}$ and
$\hat{F}$ to
\begin{subequations}
\label{eq:qLeibniz}
\begin{align}
  \hat{E} \tr (\hat{p}_1 \hat{p}_2)
  &= ( \hat{E} \tr \hat{p}_1 ) ( \E^{\h \hat{H}} \tr \hat{p}_2)
    + \hat{p}_1 (\hat{E} \tr \hat{p}_2) \\
  \hat{F} \tr (\hat{p}_1 \hat{p}_2)
  &= ( \hat{F} \tr \hat{p}_1 ) \hat{p}_2
    + (\E^{-\h \hat{H}} \tr \hat{p}_1) (\hat{F} \tr \hat{p}_2) \,,    
\end{align}
\end{subequations}
for $\hat{p}_1, \hat{p}_2 \in \Xcal_\h$. Expanding
Eqs.~\eqref{eq:qCommutation} and \eqref{eq:qLeibniz} in $\h$, we see
that in zeroth order the commutation relations and the Leibniz rule
coincide with their undeformed counterparts. This shows that the
symmetry structure of the quantum plane is a deformation of the
symmetry structure of the plane.

It is natural to ask how the deformed and undeformed symmetries are
related. This question can be answered not only for the quantum plane
but for the general case.

\subsection{The rigidity of symmetry structures}

Let $\Xcal = \mathbb{C}[x_1, \ldots, x_n]$ be the polynomial function
algebra of an $n$-dimensional space and $\Xcal_\h$ be a formal
deformation of $\Xcal$ in the sense of
Gerstenhaber \cite{Gerstenhaber:1964} with deformation parameter $\h$.
That is, $\Xcal_\h$ is an $\h$-adic algebra which is isomorphic to
$\Xcal$ modulo $\h$ as algebra, $\Xcal \cong \Xcal_\h / \h\Xcal_\h$.
In other words, setting $\h = 0$ in $\Xcal_\h$ yields $\Xcal$ as
undeformed limit. Clearly, this notion of a formal deformation as
``isomorphic modulo $\h$'' can be extended to other algebraic
structures such as a Hopf or a module structure. In addition to this
deformation property, we will assume $\Xcal_\h$ to be
$\mathbb{C}[[\h]]$-linearly isomorphic to $\Xcal[[\h]]$ as vector
space. (Mathematically speaking, this is equivalent to assuming that
$\Xcal_\h$ is topologically free.) From a computational viewpoint,
deformations without this property would have some pathological
properties. As before, we will call a $\mathbb{C}[[\h]]$-linear
isomorphism of vector spaces $\varphi: \Xcal[[\h]] \rightarrow
\Xcal_\h$ an ordering prescription.

Let us further assume that there are module structures on $\Xcal$ and
$\Xcal_\h$, denoting the module maps by
\begin{equation}
\label{eq:ModuleMaps}
  \rho: \Ug \otimes \Xcal \rightarrow \Xcal \,,\qquad
  \rho_\h: \Uhg \otimeshat \Xcal_\h \rightarrow \Xcal_\h \,,
\end{equation}
where $\mathfrak{g}$ is a semi-simple Lie algebra. We assume that the
module structure of the deformed space is a formal deformation of the
undeformed module structure, as well. This means, that the symmetry
algebra $\Uhg$ is a deformation of $\Ug$ and the action $\rho_\h$ is a
deformation of $\rho$ in the sense of ``isomorphic modulo $\h$''.  The
deformation property of $\Uhg$ turns out to be surprisingly
restrictive:
\begin{Theorem}[Gerstenhaber--Whitehead]
\label{th:RigidAlgebra}
  Let $\Ug$ be the enveloping algebra of a semi-simple Lie algebra
  $\mathfrak{g}$ and $\Uhg$ an $\h$-adic algebra which is isomorphic
  to $\Ug$ modulo $\h$, $\Uhg/\h\Uhg \cong \Ug$. Then $\Uhg$ is
  isomorphic to $\Ug[[\h]]$ as $\h$-adic algebra.
\end{Theorem}
This theorem tells us that $\Ug$ cannot be truly deformed at all.
Algebras with this property are called rigid. Gerstenhaber has
shown \cite{Gerstenhaber:1964} that an algebra is rigid if its second
Hochschild cohomology is zero. The second Whitehead lemma implies that
the second Hochschild cohomology for the enveloping algebra of a
semi-simple Lie algebra is zero, $H^2(\Ug,\Ug) = 0$.  Hence, $\Ug$ is
rigid, which was observed and used by
Drinfeld \cite{Drinfeld:1989,Drinfeld:1989b}. Note that, while the
algebra structure of $\Uhg$ is not a true deformation of $\Ug$, the
deformation of the Hopf structure of $\Uhg$ (the Leibniz rule) is
\emph{not} isomorphic to the Hopf structure of $\Ug$.

We now turn to the deformation $\rho_\h$ of the action $\rho$. Let
$\alpha: \Ug[[\h]] \rightarrow \Uhg$ be the isomorphism of $\h$-adic
algebras from theorem~\ref{th:RigidAlgebra}. Let $\varphi: \Xcal[[\h]]
\rightarrow \Xcal_\h$ be an ordering prescription. Using the
isomorphisms $\alpha$ and $\varphi$ we can proceed as for the
multiplication and transfer the deformed action $\rho_\h$ of $\Uhg$ on
$\Xcal_\h$ to an action of $\Ug[[\h]]$ on $\Xcal[[\h]]$.  Requiring
the diagram
\begin{equation}
\label{eq:StarAct}
\begin{CD}
  \Ug[[\h]] \otimeshat \Xcal[[\h]] 
  @>{\alpha \otimes \varphi }>> \Uhg \otimeshat \Xcal_\h \\
  @VV{\rho_{\alpha,\varphi}}V @VV{\rho_\h}V \\ 
  \Xcal[[\h]] @>{\varphi }>> \Xcal_\h
\end{CD}
\end{equation}
to commute, we have to define the transferred action as
\begin{equation}
\label{eq:ActionDef}
  \rho_{\alpha,\varphi} =
  \varphi^{-1} \circ \rho_\h \circ (\alpha \otimes \varphi)\,.
\end{equation}
The assumption that $\rho_\h$ is a deformation of $\rho$ in the sense
of ``isomorphic modulo $\h$'' can now be conveniently stated as
\begin{equation}
\label{eq:ActionDeform}
  \rho_{\alpha,\varphi } = \rho + \Ocal(\h) \,,
\end{equation}
which can be shown to hold independently of the choice of $\alpha$ and
$\varphi$. It turns out that the action of the symmetry algebra on
$\Xcal$ cannot be truly deformed, either:

\begin{Theorem}
\label{th:RigidAction}
Let $\rho_\h$ and $\rho$ be module maps as in
Eq.~\eqref{eq:ModuleMaps} such that $\rho_\h$ is a deformation of
$\rho$. Then there is an ordering prescription $\varphi$ such that for
$\rho_{\alpha,\varphi}$ defined as in Eq.~\eqref{eq:ActionDef} we have
$\rho_{\alpha,\varphi} = \rho$.
\end{Theorem}

\begin{proof}
  The first Whitehead lemma, which states that for semi-simple
  $\mathfrak{g}$ the first Lie algebra cohomology group vanishes,
  implies that the first Hochschild cohomology of $\Ug$ is zero, as
  well, $H^1(\Ug, \Ug) = 0$.  
  
  Consider the undeformed and deformed structure homomorphisms
  \begin{subequations}
  \label{eq:Structu}
  \begin{xalignat}{2}
    R &: \Ug \rightarrow \mathrm{End}_{\mathbb{C}[[\h]]}(\Xcal[[\h]])
    & R(g)x &:= \rho(g \otimes x) \\
    R_{\alpha,\varphi}
    &: \Ug \rightarrow \mathrm{End}_{\mathbb{C}[[\h]]}(\Xcal[[\h]])
    & R_{\alpha,\varphi}(g)x
    &:= \rho_{\alpha,\varphi}(g \otimes x)
  \end{xalignat}
  \end{subequations}
  The deformation property now reads $R_{\alpha,\varphi } = R +
  \Ocal(\h)$.  The vanishing of the first Hochschild cohomology
  implies that $R_{\alpha,\varphi }$ and $R$ are related by an inner
  automorphism (see e.g. \cite{Kassel}).  That is, there is an
  invertible $A \in \mathrm{End}_{\mathbb{C}[[\h]]}(\Xcal[[\h]])$ such
  that $R_{\alpha,\varphi}(g) = A\, R(g) \, A^{-1}$.
  
  Now we can define an ordering prescription by $\varphi' := \varphi
  \circ A$.  By definition~\eqref{eq:ActionDef} of the action we get
  $R_{\alpha,\varphi'}(g) = A^{-1} R_{\alpha,\varphi}(g)\, A = R(g)$.
  Looking at definition~\eqref{eq:Structu} of the structure maps we
  conclude that $\rho_{\alpha,\varphi'} = \rho$.
\end{proof}

Theorems~\ref{th:RigidAlgebra} and \ref{th:RigidAction} are very
general. No matter how radical a formal deformation of the symmetry
structure may seem, it is always isomorphic to the undeformed
symmetry. Most ordering prescriptions, such as the popular normal or
symmetric ordering, will only obscure this fact. But how do we find
the isomorphism of algebras of theorem~\ref{th:RigidAlgebra} and the
``good'', symmetry preserving ordering prescription of
theorem~\ref{th:RigidAction}? There is no general answer to this
question, because the elegant homological methods by which the
rigidity theorems can be proved are non-constructive.  Nevertheless,
we will demonstrate for the quantum plane and quantum Euclidean
four-space how representation theory can provide a computational
access on this problem.

\subsection{Symmetry preserving ordering prescriptions}

For each $j \in \frac{1}{2} \mathbb{N}_0$ there is an irreducible
spin-$j$ representation of the deformed symmetry algebra $\Uhsl$,
defined on the generators $\hat{E}$, $\hat{H}$, and $\hat{F}$
by \cite{Sklyanin:1985}
\begin{equation}
\label{eq:Irrepsh}
\begin{aligned}
  \hat{E}\Ket{j,m} &= \E^{\h(m+1)} \sqrt{[j+m+1][j-m]}
                           \,\Ket{j,m+1} \\
  \hat{F}\Ket{j,m} &= \E^{-\h m} \sqrt{[j+m][j-m+1]}\,\Ket{j,m-1}\\
  \hat{H}\Ket{j,m} &= 2m\Ket{j,m} \,, 
\end{aligned}
\end{equation}
on the $(2j+1)$-dimensional weight basis $\{\Ket{j,m}, m=
-j,-j+1,\ldots,j\}$. Here, each pair of brackets denotes a quantum
number,
\begin{equation}
\label{eq:qNumber}
  [a] := \frac{\E^{\h a}-\E^{-\h a}}{\E^\h-\E^{-\h}}\,,
\end{equation}
which is an $\h$-adic series in the indeterminate $a$.

The undeformed limit of Eqs.~\eqref{eq:Irrepsh} yields the spin-$j$
representation of the undeformed algebra $\Usl$ defined on the
generators by
\begin{equation}
\label{eq:Irreps}
\begin{aligned}
  E\Ket{j,m} &= \sqrt{(j+m+1)(j-m)} \,\Ket{j,m+1} \\
  F\Ket{j,m} &= \sqrt{(j+m)(j-m+1)}\,\Ket{j,m-1}\\
  H\Ket{j,m} &= 2m\Ket{j,m} \,. 
\end{aligned}
\end{equation}
Let us formally define the operators $M := \tfrac{1}{2}H$ and
\begin{equation}
  J := \tfrac{1}{2}(\sqrt{2C + 1} - 1)\,,
  \quad\text{where}\quad
  C := EF + FE + \tfrac{1}{2}H^2 \,,
\end{equation}
such that $M\Ket{j,m} = m\Ket{j,m}$ and $J\Ket{j,m} = j\Ket{j,m}$.
Remark, that $J$ is not an element of $\Ug$ proper. Now we can define
\begin{equation}
\label{eq:alphaDef}
\begin{aligned}
  \alpha^{-1}(\hat{E}) &:= E \E^{\h(M+1)}
  \sqrt{\frac{[J+M+1][J-M]}{(J+M+1)(J-M)}}\\
  \alpha^{-1}(\hat{F}) &:= F \E^{-\h M}
  \sqrt{\frac{[J+M][J-M+1]}{(J+M)(J-M+1)}}\\
  \alpha^{-1}(\hat{H}) &:= H 
\end{aligned}
\end{equation}
where the right hand sides have to be understood as $\h$-adic series
with polynomials in $J$ and $M$ as coefficients. Since all expressions
involved are symmetric with respect to $J \mapsto -J-1$, the operator
$J$ appears in the coefficient polynomials only as polynomial of
$2J(J+1) = C$. For example,
\begin{equation}
  \alpha^{-1}(\hat{E}) = E\bigl\{ 1 + \tfrac{1}{2}( 2 + H ) \h +
  \tfrac{1}{12}[ C + (1 + H)( 5 + 2H)] \h^2 \bigr\} + \Ocal(\h^3) \,. 
\end{equation}
We conclude, that the operators defined in Eq.~\eqref{eq:alphaDef} can
be viewed as elements of $\Usl[[\h]]$.

Comparing Eqs.~\eqref{eq:Irrepsh} and \eqref{eq:Irreps} we see that
the operators $\alpha^{-1}(\hat{E})$, $\alpha^{-1}(\hat{F})$, and
$\alpha^{-1}(\hat{H})$ have the same irreducible representations as
$\hat{E}$, $\hat{F}$, and $\hat{H}$, respectively.  Hence,
Eqs.~\eqref{eq:alphaDef} define a homomorphism of algebras
$\alpha^{-1}: \Uhsl \rightarrow \Usl[[\h]]$. Analogously to
Eqs.~\eqref{eq:alphaDef} we can construct the inverse of
$\alpha^{-1}$, which shows that $\alpha^{-1}$ is an isomorphism of
algebras. Its inverse $\alpha$ is the searched-for isomorphism of
theorem~\ref{th:RigidAlgebra}.

In order to find the symmetry preserving ordering prescription of
theorem~\ref{th:RigidAction} we want to decompose the plane $\Xcal$
and the quantum plane $\Xcal_\h$, as defined in
Eq.~\eqref{eq:PlanesDef}, into irreducible subrepresentations with
respect to the $\Usl[[\h]] \cong \Uhsl$ module structure. Define the
representation matrices by
\begin{equation}
  R^j(g)^{m}{}_{m'} := \Bra{j,m} g \Ket{j,m'} \quad\text{and}\quad
  R^j_\h(\hat{g})^{m}{}_{m'} := \Bra{j,m} \hat{g} \Ket{j,m'} 
\end{equation}
for all $g \in \Usl$, $\hat{g} \in \Uhsl$. Note, that the isomorphism
$\alpha: \Usl[[\h]] \rightarrow \Usl$ was defined in
Eqs.~\eqref{eq:alphaDef} precisely such that
\begin{equation}
\label{eq:Requal}
  R^j(g)^{m}{}_{m'} = R^j_\h(\alpha(g))^{m}{}_{m'} \,.
\end{equation}
The basis of $\Xcal_\h$ which decomposes the quantum plane into
irreducible subrepresentations can be computed to \cite{Blohmann:2002a}
\begin{equation}
\label{eq:hTensorBasis}
  \hat{T}^j_m :=
  \qBinom{2j}{j+m}^{\frac{1}{2}} \hat{x}^{j-m} \hat{y}^{j+m}
  \,,\quad\text{where}\quad
  \qBinom{n}{k} :=
   \E^{\h k(k-n)}\,\frac{[n]!}{[n-k]![k]!}
\end{equation}
is a deformation of the binomial coefficient (sticking to the standard
notation of \cite{Schmuedgen}) and $[n]! := [1][2]\cdots[n]$
for natural $n$.  The expression for the basis which decomposes the
commutative plane is the undeformed limit
\begin{equation}
%\label{eq:hTensorBasis}
  T^j_m :=
  \tbinom{2j}{j+m}^{\frac{1}{2}} x^{j-m} y^{j+m} \,.
\end{equation}
The action of the symmetry algebras on these bases is 
\begin{equation}
\label{eq:IrredAction}
  \rho_\h(\hat{g} \otimes \hat{T}^j_m)
  = \hat{T}^j_{m'} R_\h^j(\hat{g})^{m'}{}_m
  \quad\text{and}\quad
  \rho(g \otimes T^j_m) = T^j_{m'} R^j(g)^{m'}{}_m \,,
\end{equation}
where we sum over repeated indices. Defining an ordering prescription
by
\begin{equation}
\label{eq:Ordering}
  \varphi(T^j_{m}) = \hat{T}^j_{m} \quad\Leftrightarrow\quad
  \varphi(x^k y^l) =
  \qBinom{k+l}{k}^{\frac{1}{2}}
  \tbinom{k+l}{k}^{-\frac{1}{2}}
  \hat{x}^k \hat{y}^l \,,
\end{equation}
we find that for all $g \in \Ug$ we have
\begin{equation}
\begin{split}
  \rho(g \otimes T^j_m)
  &= T^j_{m'} R^j(g)^{m'}{}_m 
  = (\varphi^{-1}\circ \varphi)
    \bigl(T^j_{m'} R^j(g)^{m'}{}_m \bigr) \\
  &= \varphi^{-1}
     \bigl(\hat{T}^j_{m'} R_\h^j(\alpha(g))^{m'}{}_m \bigr)
  = \varphi^{-1}
    \bigl( \rho_\h(\alpha(g) \otimes \hat{T}^j_{m'}) \bigr) \\
  &= (\varphi^{-1} \circ \rho_\h \circ [\alpha \otimes \varphi])
    (g \otimes T^j_{m'} ) 
  = \rho_{\alpha,\varphi}(g \otimes T^j_{m'}) \,,
\end{split}
\end{equation}
where we used Eqs.~\eqref{eq:IrredAction}, \eqref{eq:Ordering},
\eqref{eq:Requal}, and \eqref{eq:ActionDef}. Since $\{ T^j_m\}$ is a
basis of $\Xcal$, it follows that $\rho_{\alpha,\varphi} = \rho$.
Hence, the ordering prescription~\eqref{eq:Ordering} is the
searched-for symmetry preserving ordering of
theorem~\ref{th:RigidAction}.

\subsection{Quantum Euclidean four-space}

Finally, we want to give the result of the analogous computations for
quantum Euclidean four-space.

The algebra of a commutative four-dimensional space is the polynomial
algebra generated by its four coordinates. For convenience, the
coordinates can be arranged in a matrix $\bigl(\begin{smallmatrix} a &
  b \\ c & d
  \end{smallmatrix}\bigr) :=
  \bigl(\begin{smallmatrix}
    x_0 - \I x_3 & \I x_1 +  x_2 \\ \I x_1 -  x_2 & x_0 + \I x_3
\end{smallmatrix}\bigr)$
such that square of the invariant four-length $l^2 := x_0^2 + x_1^2 +
x_2^2 + x_3^2$ is given by the determinant. The Euclidean four-space
algebra can now be viewed as the algebra of 2$\times$2-matrices
matrices $M(2) := \mathbb{C}[a,b,c,d]$.

Quantum Euclidean four-space is defined as the $\h$-adic algebra
$M_\h(2)$ of quantum 2$\times$2-matrices, generated by $\hat{a}$,
$\hat{b}$, $\hat{c}$, $\hat{d}$ with relations
\begin{equation}
\label{eq:MqRel}
\begin{gathered}
  \hat{a}\hat{b} = \E^\h \hat{b}\hat{a}\,,
  \quad \hat{a}\hat{c} = \E^\h \hat{c} \hat{a}\,,
  \quad \hat{b}\hat{d} = \E^\h \hat{d}\hat{b}\,,
  \quad \hat{c}\hat{d} = \E^\h \hat{d}\hat{c} \\
  \hat{b}\hat{c} = \hat{c}\hat{b}\,,\quad
  \hat{a}\hat{d} - \hat{d}\hat{a} = (\E^\h - \E^{-\h}) \hat{b}\hat{c} \,.
\end{gathered}
\end{equation}
The commutation relations of the usual quantum matrices for a real
deformation parameter \cite{Woronowicz:1987b} can be obtained by the
substitution $\E^\h \mapsto q$.  The central and invariant square
$\hat{l}^2$ of the quantum four-length is given by the quantum
determinant
\begin{equation}
  \hat{l}^2 := \hat{a}\hat{d} - \E^{\h} \hat{b}\hat{c} \,.
\end{equation}
Quantum Euclidean four-space carries by construction a representation
of the quantum orthogonal algebra $\Uhso$, which is, analogous to the
undeformed case, the tensor algebra of two copies of $\Uhsl$, $\Uhso =
\Uhsl \otimeshat \Uhsl$. This shows that the isomorphism $\beta$ from
the orthogonal algebra to the quantum orthogonal algebra of
theorem~\ref{th:RigidAlgebra} is simply given by the tensor square of
the isomorphism $\alpha: \Usl[[\h]] \rightarrow \Uhsl$, which was
constructed in the last section,
\begin{equation}
  \beta: \Uso[[\h]] \stackrel{\cong}{\longrightarrow} \Uhso
  \,,\quad\text{where}\quad \beta := \alpha \otimes \alpha \,.
\end{equation}

In order to calculate the symmetry preserving ordering prescription we
again need to reduce the deformed and the undeformed Euclidean space
into their irreducible subrepresentations. This time the space
algebras each possess a non-trivial invariant element, $l^2$ and
$\hat{l}^2$, so the irreducible subspaces are degenerate. More
precisely, every highest weight vector of $M_\h(2)$ is of the form
$\hat{z}\hat{d}^{2j}$ for $2j \in \mathbb{N}_0$, where $\hat{z} \in
\mathbb{C}[\hat{l}^2] \subset M_\h(2)$, the weight being $(j,j)$. Let
us denote by $\{\hat{T}^{(j,j)}_{mm'}\}$ the basis of the irreducible
$(j,j)$-subrepresentation of $M_\h(2)$ which is generated by
$\hat{d}^{2j}$ and by $\{T^{(j,j)}_{mm'}\}$ the according basis of
$M(2)$ generated by $d^{2j}$.

If we want the symmetry preserving ordering prescription $\varphi$
additionally to preserve the degree of the monomials, then $\varphi$
must identify these bases, $\varphi(T^{(j,j)}_{mm'}) =
\hat{T}^{(j,j)}_{mm'}$. Moreover, as $l^2$ and $\hat{l}^2$ are the
only invariant elements of degree 2, we must have $\varphi(l^2) \sim
\hat{l}^2$.  For convenience we choose the proportionality constant to
be 1. Observing that since $l^2$ is invariant we have $\varphi(l^{2n}
x) = \varphi(l^{2n}) \varphi(x) = \hat{l}^{2n} \varphi(x)$ for all $x
\in M(2)$, we conclude that the symmetry preserving ordering
prescription must be defined as
\begin{equation}
\label{eq:phi0}
  \varphi(l^{2n} T^{(j,j)}_{mm'})
  = \hat{l}^{2n}\hat{T}^{(j,j)}_{mm'} \,.
\end{equation}

We now want to express this ordering prescription in terms of the
normal ordered (Poincar\'e--Birkhoff--Witt) bases. Towards this end we
need to expand the irreducible basis $\hat{T}^{(j,j)}_{mm'}$ in terms
of the normal ordered basis, and vice versa.  Our starting point will
be the multiplication map of $M_\h(2)$ in terms of the irreducible
basis, which has been calculated explicitly in \cite{Blohmann:2002a},
\begin{equation}
\label{eq:EuclProd1}
  \hat{T}^{(j_1,j_1)}_{m_1m_1'} \hat{T}^{(j_2,j_2)}_{m_2m_2'} 
  = \sum_{j,m,m'} \CGqs{j_1}{j_2}{j}{m_1}{m_2}{m}
  \CGqs{j_1}{j_2}{j}{m_1'}{m_2'}{m'} \,
  \hat{l}^{2(j_1+j_2-j)} \, \hat{T}^{(j,j)}_{m, m'} \,,
\end{equation}
where the expressions in parentheses denote the $\h$-adic quantum
Clebsch--Gordan coefficients of $\Uhsl$, which can be obtained from
the $q$-deformed Clebsch--Gordan coefficient \cite{Schmuedgen} by the
substitution $q \mapsto \E^\h$. Observing, that
\begin{equation}
\label{eq:Textreme}
  \hat{T}^{(j,j)}_{-j, -j} = \hat{a}^{2j} \,,\quad
  \hat{T}^{(j,j)}_{-j, j} = \hat{b}^{2j} \,,\quad
  \hat{T}^{(j,j)}_{j, -j} = \hat{c}^{2j} \,,\quad
  \hat{T}^{(j,j)}_{j, j} = \hat{d}^{2j}
\end{equation}
we get from Eq.~\eqref{eq:EuclProd1}
\begin{equation}
\label{eq:Textreme2}
\begin{aligned}
  \hat{a}^{2n_a} \hat{b}^{2n_b} 
  &= \qBinom{2n_a + 2n_b}{2n_b}^{-\frac{1}{2}}
  \hat{T}^{(n_a+n_b,n_a+n_b)}_{-n_a-n_b,-n_a+n_b}\\
  \hat{c}^{2n_c} \hat{d}^{2n_d} 
  &= \qBinom{2n_c + 2n_d}{2n_c}^{-\frac{1}{2}}
  \hat{T}^{(n_c+n_d,n_c+n_c)}_{n_c+n_d,-n_c+n_d}\,.
\end{aligned}
\end{equation}
Multiplying these two monomials by Eq.~\eqref{eq:EuclProd1} we get an
expression of the normal ordered basis in terms of the reduced basis,
\begin{equation}
\begin{split}
\label{eq:phi1}
  \hat{a}^{2n_a} \hat{b}^{2n_b} \hat{c}^{2n_c} \hat{d}^{2n_d}
  &=\sum_j
    \CGqs{n_a+n_b}{n_c+n_d}{j}{-n_a-n_b\,}{n_c + n_d }{n_1} 
    \CGqs{n_a+n_b}{n_c+n_d}{j}{-n_a+n_b\,}{-n_c + n_d }{n_2}\\
  &\qquad\times
    \qBinom{2n_a + 2n_b}{2n_b}^{-\frac{1}{2}}
    \qBinom{2n_c + 2n_d}{2n_c}^{-\frac{1}{2}}
    \hat{l}^{2(n - j)} \, \hat{T}^{(j,j)}_{n_1 n_2}
\end{split}
\end{equation}
for all $n_a, n_b, n_c, n_d \in \frac{1}{2}\mathbb{N}_0$, where $n_1
:= -n_a - n_b + n_c + n_d$, $n_2 := -n_a + n_b - n_c + n_d$, and $n :=
n_a + n_b + n_c + n_d$. Using the orthogonality of the quantum
Clebsch--Gordan coefficients Eq.~\eqref{eq:phi1} can be inverted,
\begin{equation}
\begin{split}
\label{eq:phi2}
    \hat{l}^{2(n - j)} \, \hat{T}^{(j,j)}_{n_1 n_2}
  &=\sum_{k}
    \CGqs{n_a+n_b}{n_c+n_d}{j}{-n_a-n_b\,}{n_c + n_d }{n_1}^{-1} 
    \CGqs{n_a+n_b}{n_c+n_d}{j}{-n_a+n_b+k\,}{-n_c+n_d-k}{n_2}\\
  &\qquad\times
    \qBinom{2n_a + 2n_b}{2n_b+k}^{\frac{1}{2}}
    \qBinom{2n_c + 2n_d}{2n_c+k}^{\frac{1}{2}}
    \hat{a}^{2n_a-k} \hat{b}^{2n_b+k} \hat{c}^{2n_c+k} \hat{d}^{2n_d-k} \,.
\end{split}
\end{equation}
Applying the ordering prescription~\eqref{eq:phi0} to the undeformed
limit of Eq.~\eqref{eq:phi1} and inserting Eq.~\eqref{eq:phi2} we
finally obtain
\begin{equation}
\begin{split}
\label{eq:phi3}
  \varphi(a^{2n_a} b^{2n_b} c^{2n_c} d^{2n_d})
  &=\sum_{j,k}
    \frac{
      \CGs{n_a+n_b}{n_c+n_d}{j}{-n_a-n_b\,}{n_c + n_d }{n_1}
      \qBinom{2n_a + 2n_b}{2n_b}^{\frac{1}{2}}
      \qBinom{2n_c + 2n_d}{2n_c}^{\frac{1}{2}}
    }{%
      \CGqs{n_a+n_b}{n_c+n_d}{j}{-n_a-n_b\,}{n_c + n_d }{n_1}
      \tbinom{2n_a + 2n_b}{2n_b+k}^{\frac{1}{2}}
      \tbinom{2n_c + 2n_d}{2n_c+k}^{\frac{1}{2}}
    } \\
  &\qquad\times  
    \CGs{n_a+n_b}{n_c+n_d}{j}{-n_a+n_b\,}{-n_c + n_d }{n_2}
    \CGqs{n_a+n_b}{n_c+n_d}{j}{-n_a+n_b+k\,}{-n_c+n_d-k}{n_2}\\
  &\qquad\times
   \hat{a}^{2n_a-k} \hat{b}^{2n_b+k} \hat{c}^{2n_c+k} \hat{d}^{2n_d-k} \,,
\end{split}
\end{equation}
which is the desired expression for the symmetry preserving ordering
prescription in terms of the normal ordered bases.

Via the $q$-binomial and $q$-Clebsch--Gordan coefficients, the
dependence on the deformation parameter $\h$ is contained in basic
hypergeometric series. In principle, the expansion of
Eq.~\eqref{eq:phi3} in powers of $\h$ could be expressed in terms of
combinatorial partition functions. Since there are no general number
theoretic formulas for partition functions, this would not be of much
practical value. Trying to calculate the coefficients of the
$\h$-expansion explicitly, one would quickly learn that basic
hypergeometric functions are the more efficient way to deal with
partitions. Therefore, Eq.~\eqref{eq:phi3} is probably already the
form which is best suited for applications. Moreover, using
Eq.~\eqref{eq:phi3} as generating functional, the $\h$-expansion can
be done by computer algebra.

\section{Conclusion}

In this paper perturbative deformations of symmetry structures on
noncommutative spaces were studied. It was shown that the rigidity of
symmetry algebras extends to rigidity of the action of the symmetry on
the space. This result applies to all spaces with a symmetry given by
a semi-simple Lie algebra and, hence, comprises most spaces which are
of interest in physics. The generality of the results may be
surprising at first sight: Even if a formal deformation of the
symmetry structure looks extremely complicated, it is always
isomorphic to the undeformed symmetry. But one has to keep in mind
that the class of isomorphisms between rings over formal power series
is very large. In general, the isomorphisms between the deformed and
undeformed symmetry structure will not make numerical sense for a
particular value of the deformation parameter. For the physical
interpretation, this is not necessarily a problem as long as one
stays within the realm of perturbation theory. This situation is not
much worse than for ordinary quantum field theory where true
convergence of loop expansions cannot be obtained easily, if at all.

The results obtained here can be applied to the construction of gauge
theories on noncommutative spaces \cite{Madore:2000b,Jurco:2001} which
is situated entirely in the realm of formal power series. In this
context, the rigidity of symmetry structures has interesting
implications for the construction of invariants, which would have to
appear in Lagrangians. Consider the deformation of a space with a
deformed symmetry structure such as quantum Euclidean space. If the
star product is implemented with the symmetry preserving ordering
prescription, then invariance with respect to the deformed symmetry is
the same as invariance with respect to the undeformed symmetry, since
the symmetry preserving ordering prescription maps invariants to
invariants. Moreover, invariants can be constructed using the quantum
metric.

This may not look surprising but let us illustrate with an example
what can go wrong. Since the quantum plane does not have nontrivial
invariants we take quantum Euclidean four-space. If the star product
were realized by normal ordering the quantum Euclidean four-length
would be given by
\begin{equation}
  a \star d - \E^\h b \star c
  = \varphi_{\mathrm{normal}}^{-1}
  (\hat{a}\hat{d} - \E^\h \hat{b} \hat{c} )
  = ad - \E^\h bc 
\end{equation}
which is not the invariant $ad-cd = x_0^2 + x_1^2 + x_2^2 + x_3^2$.
In contrast, the ordering prescription given by formula
Eq.~\eqref{eq:phi3} yields $\varphi(ad) = \hat{a}\hat{b}$ and
$\varphi(bc) = \E^\h \hat{b} \hat{c}$. Hence $\varphi^{-1}
(\hat{a}\hat{d} - \E^\h \hat{b} \hat{c} ) = ad - bc$, as claimed.
While it would be simple to ad-hoc modify the normal ordering in such
a way that the quadratic invariant is preserved, preserving invariants
of all orders would be difficult.

\providecommand{\href}[2]{#2}\begingroup\raggedright\endgroup

\end{document}